\documentclass[12pt]{article}

\usepackage{amsmath}
\usepackage{amssymb}
\usepackage{eucal}

\usepackage[english]{babel}

\begin{document}

\baselineskip 16pt

\title{ON MAXIMAL SUBGROUPS OF A FINITE SOLVABLE GROUP}

%\medskip

\author{D.\,V. Gritsuk and V.\,S. Monakhov}

%\date{}

\maketitle

%\medskip

\begin{abstract}
The following result is received:
Let $H$ be a non-normal maximal subgroup of a
finite solvable group $G$ and let $q \in \pi(F(H/\mathrm{Core}
_GH))$, then $G$ has a Sylow $q$\nobreakdash-\hspace{0pt}subgroup
$Q$ such that ${N_{G}(Q) \subseteq H}$. 
\end{abstract}


MSC2010 20D20, 20E34

%\bigskip
\section{Introduction}

All groups considered in this paper will be finite. All notations
and definitions correspond [1].

In 1986, V.\,A. Vedernikov obtained the following result.

\medskip

 { \bf Theorem A.} [2, Corollary 2.1] {\sl If $H$ is a
non-normal maximal subgroup of a solvable group $G$ then
$N_{G}(Q) \subseteq H$ for some Sylow subgroup $Q$ of $G$.}

\medskip

Here $N_{G}(Q)$ is the normalizer of $Q$ in $G$.

In this note we consider the following question:

\medskip

{\sl What is Sylow subgroup that its normalizer contained in
non-normal maximal subgroup of solvable group?}

\medskip

The answer to this question was obtained in the following theorem.

\medskip

{ \bf Theorem 1.} {\sl Let $H$ be a non-normal maximal subgroup of a
solvable group $G$ and let $q \in \pi(F(H/\mathrm{Core} _GH))$.
Then $G$ has a Sylow $q$\nobreakdash-\hspace{0pt}subgroup $Q$ such
that ${N_{G}(Q) \subseteq H}$.}

\medskip

Here $F(X)$ is the Fitting subgroup of $X$, $\pi (Y)$ is the set of
all prime divisors of $|Y|$, $\mbox{Core} _GH=\cap_{g\in G}H^g$ is 
the core of $H$ in $G$, i.\,e. the largest normal subgroup of
$G$ contained in~$H$.

\medskip

{\bf Corollary 1.1.} {\sl Let $H$ be a non-normal maximal subgroup
of a solvable group $G$ and let $q \in \pi(F(H/\emph{Core} _GH))$.
Then $H$ has a Sylow $q$\nobreakdash-\hspace{0pt}subgroup $Q$ such
that $N_{G}(H_1)\subseteq H$ for each subgroup $H_1$ of $H$
satisfying $Q\subseteq H_1\subseteq H$.}

\medskip

{\bf Corollary 1.2.} {\sl Let $H$ be a non-normal maximal subgroup
of a solvable group $G$ and let $\omega \subset \pi (F(H/\mathrm{Core}
_GH))$. Then $G$ has a Hall $\omega$\nobreakdash-\hspace{0pt}subgroup 
$G_{\omega}$ such that ${N_{G}(G_{\omega}) \subseteq H}$.}
\medskip

For non-solvable groups, this result is false. For example,
$PSL(2,17)$ has the order $2^4 \cdot 3^2 \cdot 17$ and the symmetric
group $S_{4}$ is a maximal subgroup in $PSL(2,17)$, see [2]. Since
$|S_4|=2^3\cdot 3$, it follows that $S_4$ does not contain a Sylow
subgroup of $PSL(2,17)$. Thus it is not possible to extend theorem
of V.\,A. Vedernikov and Theorem 1 to the non-solvable groups.

The following question is contained in [2].

\medskip

{\bf Question.} (V.\,A. Vedernikov, 1986, [2]) {\sl Can Theorem A be 
extended to a $p$-solvable group $G$ containing a maximal subgroup $M$ 
such that $|G:M|=p^a$?}

\medskip

The positive answer to this question was obtained for
$F(M/\mathrm{Core}_{G}M) \ne 1$ in the following theorem.

\medskip

{\bf Theorem 2.} {\sl Let $G$ be a
$p$\nobreakdash-\hspace{0pt}solvable group. Let $M$ be a non-normal
maximal subgroup of $G$, and let $|G:M|=p^{\alpha}$. Then:

$1)$ if $F(M/\mathrm{Core}_{G}M) \ne 1$ and $q \in
\pi(F(M/\mathrm{Core}_{G}M))$, then $G$ has a Sylow
$q$\nobreakdash-\hspace{0pt}subgroup $Q$ such that $N_{G}(Q)
\subseteq M$;

$2)$ if $F(M/\mathrm{Core}_{G}M)=1$, then $N_{G}(K) \subseteq M$ for
some Hall $p^{\prime}$\nobreakdash-\hspace{0pt}subgroup $K$ of $G$.}

%\medskip

\section{Notations and preliminary results}

%\medskip

In this section we give some definitions and basic results that will
be used later in our paper.

Let $\Bbb P$ be the set of all prime numbers, and let $\pi$ be the set
of primes, i.\,e $\pi \subseteq \Bbb P$. In that paper,
$\pi^{\prime}$ is the set of all primes not contained in $\pi$,
i.\,e. $\pi= \Bbb P \setminus \pi^{\prime}$, $\pi(m)$ is the set of
prime divisors of $m$. If $\pi(m) \subseteq \pi$ then $m$ is called
{\sl $\pi$\nobreakdash-\hspace{0pt}number}.

The subgroup $H$ of $G$ is called {\sl
$\pi$\nobreakdash-\hspace{0pt}subgroup} if $|H|$ is
$\pi$\nobreakdash-\hspace{0pt}number. The subgroup $H$ of $G$ is
called Hall {\sl $\pi$\nobreakdash-\hspace{0pt}subgroup} if $|H|$ is
$\pi$\nobreakdash-\hspace{0pt}number and $|G:H|$ is
$\pi^{\prime}$\nobreakdash-\hspace{0pt}number. As usually,
$O_{\pi}(X)$ is the largest normal
$\pi$\nobreakdash-\hspace{0pt}subgroup of $X$.

The group is called {\sl $\pi$\nobreakdash-\hspace{0pt}separable} if
it has a series whose factors are
$\pi$\nobreakdash-\hspace{0pt}groups or
$\pi^{\prime}$\nobreakdash-\hspace{0pt}groups.

{\bf Lemma 1.} [2, Theorem 1] {\sl Let $G$ be a $\pi$-separable
group, and let $H$ be a subgroup of $G$. If $|G:H|$ is $\pi$-number,
then $O_{\pi}(H)\subseteq O_{\pi}(G)$.}

{\bf Lemma 2.} {\sl Let $R$ be a Hall $\pi$-subgroup of a
$\pi$-separable group $G$, and let $N$ be a normal subgroup of $G$.
Then $N_{G}(R)N/N=N_{G/N}(RN/N)$.}

%\medskip

{\sc Proof.} For $x \in N_{G}(R)$ we obtain:
$$
(x^{-1}N)(RN/N)(xN)=R^{x}N/N=RN/N,
$$
i.\,e. $N_{G}(R)N/N \leq N_{G/N}(RN/N)$. Conversely, if $yN \in
N_{G/N}(RN/N)$, then $R^{y}N=RN$. Now $R$ and $R^{y}$ are Hall
subgroups of $RN$ which are conjugate, i.\,e.
$R^{y}=R^{ak}=R^{k}$ for some $ak \in RN$, $a \in R$, $k \in N$.
Then $yk^{-1} \in N_{G}(R)$, whence $y \in N_{G}(R)N$, i.\,e.
$N_{G/N}(RN/N) \leq N_{G}(R)N/N$. The Lemma is proved.

{\bf Lemma 3.} {\sl Let $G$ be a
$\pi$\nobreakdash-\hspace{0pt}solvable group containing a nilpotent
Hall $\pi$\nobreakdash-\hspace{0pt}subgroup. If $H$ is a maximal
subgroup of $G$ and $|G:H|$ is $\pi$\nobreakdash-\hspace{0pt}number,
then $O_{\pi}(H)$ is a normal subgroup of $G$.}

%\medskip

{\sc Proof.} By Lemma 1, $O_{\pi}(H) \subseteq O_{\pi}(G)$. If
$O_{\pi}(H)=O_{\pi}(G)$ then $O_{\pi}(H)$ is a normal subgroup of
$G$. Let $O_{\pi}(H)$ be a proper subgroup of $O_{\pi}(G)$, and let
$G_{\pi}$ be a Hall $\pi$\nobreakdash-\hspace{0pt}subgroup of $G$.
Clearly, $O_{\pi}(G)$ is a proper subgroup of $G_{\pi}$. Since
$G_{\pi}$ is a nilpotent subgroup, we have $O_{\pi}(H)$ is a
proper subgroup of $D=N_{G_{\pi}}(O_{\pi}(H))$. Since $O_{\pi}(H)=H
\cap O_{\pi}(G)$, we have $D$ does not contained in $H$, so
$N_{G}(O_{\pi}(H))\supseteq \langle H,D\rangle =G$ and $O_{\pi}(H)$
is a normal subgroup of $G$. The Lemma is proved.

%Note that conclusion 1 of a Corollary 1.2 [2] follows from Lemma 3
%if $\pi = \{p\}$.

%\bigskip

%\medskip
\section{Main results}

{\bf Theorem 3.} {\sl Let $G$ be a $\pi$-solvable group containing 
a nilpotent Hall $\pi$-subgroup. Let $M$ be a non-normal maximal subgroup of $G$, 
and let $|G:M|$ is $\pi$\nobreakdash-\hspace{0pt}number. Then:

$1)$ if $F(M/\emph{Core}_{G}M) \ne 1$ and $q \in
\pi(F(M/\emph{Core}_{G}M))$, then $G$ has a Sylow
$q$\nobreakdash-\hspace{0pt}subgroup $Q$ such that $N_{G}(Q)
\subseteq M$;

$2)$ if $F(M/\emph{Core}_{G}M)= 1$, then $N_{G}(K) \subseteq M$ for
some Hall $\pi^{\prime}$\nobreakdash-\hspace{0pt}subgroup $K$ of
$G$.}

\medskip

{\sc Proof.} Case 1: $\mbox{Core}_{G}M=1$. Since $M$ is a
$\pi$\nobreakdash-\hspace{0pt}solvable group, we have  $M$
contains Hall $\pi^{\prime}$\nobreakdash-\hspace{0pt}subgroup of
$K$. By Theorem 1.39 [1], we have
$$
|G:K|=|G:M||M:K|.
$$
Since $|G:M|$ is a $\pi$\nobreakdash-\hspace{0pt}number, we have
$K$ is a Hall $\pi^{\prime}$\nobreakdash-\hspace{0pt}subgroup
of $G$. Now we have $O_{\pi^{\prime}}(G) \leq K \leq M$, so
$O_{\pi^{\prime}}(G) \leq \mbox{Core}_{G}M=1$. Since $G$ is a
$\pi$\nobreakdash-\hspace{0pt}solvable group, we have 
$O_{\pi}(G) \ne 1$. But $G$ is a primitive group with a maximal
subgroup $M$ such that $\mbox{Core}_{G}M=1$. Therefore, for some $p
\in \pi(G)$
$$
 N=O_{p}(G)=F(G)=C_{G}(O_{p}(G))\ne 1, \ G=M[O_{p}(G)], \  \Phi(G)=1
$$
and $O_{p}(M)=1$.

%\medskip

1. Let $F(M) \ne 1$, $q \in \pi(F(M))$, and let $Q$ be a Sylow
$q$\nobreakdash-\hspace{0pt}subgroup of $M$. Then $O_{q}(M) \ne 1$,
$O_{q}(M) \subseteq Q$, and since $\mbox{Core}_{G}M=1$, we have 
$N_{G}(O_{q}(M))=M$. Since $O_{p}(M)=1$, it follows that $p$ does not
belong to $\pi(F(M))$ and $p \ne q$. The subgroup $D=N_{G}(Q) \cap
O_{p}(G)$ is a normal subgroup of $N_{G}(Q)$ and
$$
D \subseteq C_{G}(Q)\subseteq C_{G}(O_{q}(M)) \subseteq
N_{G}(O_{q}(M))=M.
$$
Now $D=N_{G}(Q) \cap O_{p}(G) \subseteq M \cap O_{p}(G)=1$. We
consider the subgroup $L=N_{G}(Q)O_{p}(G)$. Since $N_{G}(Q) \cap
O_{p}(G)=1$, we have $L=[O_{p}(G)]N_{G}(Q)$. It follows from
$G=[O_{p}(G)]M$ and Dedekind identity that
$$
L=[O_{p}(G)]N_{G}(Q)=[O_{p}(G)](L \cap M), \ L/O_{p}(G) \simeq
N_{G}(Q) \simeq L \cap M,
$$
and $L\cap M$ is a $q$\nobreakdash-\hspace{0pt}closed subgroup. But
$Q \subseteq L \cap M$, so $Q$ is a normal subgroup of $L \cap M$.
It follows from $|N_{G}(Q)|=|L \cap M|$ that $L \cap M=N_{G}(Q)$ and
$N_{G}(Q) \subseteq M$.

2. Let $F(M)=1$, and let $K$ be a Hall
$\pi^{\prime}$\nobreakdash-\hspace{0pt}subgroup of $M$. By Lemma 2,
$O_{\pi}(M)$ is a normal subgroup of $G$. Hence $O_{\pi}(M) \leq
\mbox{Core}_{G}M=1$. Since $G$ is a
$\pi$\nobreakdash-\hspace{0pt}solvable group, we have
$O_{\pi^{\prime}}(M) \ne 1$ and $O_{\pi^{\prime}}(M) \subseteq K$.
Since $\mbox{Core}_{G}M=1$, we have 
$N_{G}(O_{\pi^{\prime}}(M))=M$. The subgroup $B=N_{G}(K) \cap
O_{p}(G)$ is a normal subgroup of $N_{G}(K)$ and
$$
B \subseteq C_{G}(K)\subseteq C_{G}(O_{\pi^{\prime}}(M))\subseteq
N_{G}(O_{\pi^{\prime}}(M))=M.
$$
Now $B=N_{G}(K) \cap O_{p}(G) \subseteq M \cap O_{p}(G)=1$. We
consider the subgroup $T=N_{G}(K)O_{p}(G)$. Since $N_{G}(K) \cap
O_{p}(G)=1$, we have $T=[O_{p}(G)]N_{G}(K)$. It follows from
$G=[O_{p}(G)]M$ and Dedekind identity that
$$
T=[O_{p}(G)]N_{G}(K)=[O_{p}(G)](T \cap M), \ T/O_{p}(G) \simeq
N_{G}(K) \simeq T \cap M
$$
and $T \cap M$ is a $\pi^{\prime}$\nobreakdash-\hspace{0pt}closed
subgroup. But Hall $\pi^{\prime}$\nobreakdash-\hspace{0pt}subgroup
$K$ is contained in $T \cap M$, so $K$ is a normal subgroup of $T
\cap M$ and $T \cap M \subseteq N_{G}(K)$. It follows from
isomorphism $N_{G}(K) \simeq T \cap M$ that $|N_{G}(K)|=|T \cap M|$,
so $T \cap M= N_{G}(K)$ and $N_{G}(K) \subseteq M$. The Case 1 is
proved.

Case 2:  $N= \mbox{Core}_{G}M \ne 1$. We consider the quotient group
$\overline{G}=G/N$. Clearly, $\overline{G}$ is a
$\pi$\nobreakdash-\hspace{0pt}solvable group and $|\overline{G}:
\overline{M}|=|G:M|$ is a $\pi$\nobreakdash-\hspace{0pt}number,
where $\overline{M}=M/N$. Since
$\mbox{Core}_{\overline{G}}\overline{M}=1$, for the group
$\overline{G}$ with a non-normal maximal subgroup $\overline{M}$, we
apply Case 1.

1. Let $F(\overline{M}) \ne 1$, and let $q \in
\pi(F(\overline{M}))$. By Case 1, $\overline{G}$ has a Sylow
$q$\nobreakdash-\hspace{0pt}subgroup $\overline{Q}$ such that
$N_{\overline{G}}(\overline{Q}) \subseteq \overline{M}$. Let
$\overline{Q}=A/N$, and let $Q$ be a Sylow
$q$\nobreakdash-\hspace{0pt}subgroup of $A$. Then $QN/N=A/N=
\overline{Q}$ and by Theorem 1.65 [1]
$$
N_{\overline{G}}(\overline{Q})=N_{G/N}(QN/N)=N_{G}(Q)N/N.
$$
By the condition $N_{G}(Q)N/N \subseteq \overline{M}=M/N$, so
$N_{G}(Q) \subseteq M$.

2. Let $F(M/\mbox{Core}_{G}M)=1$. By Case 1,
$N_{\overline{G}}(\overline{K}) \subseteq \overline{M}$ for some
Hall $\pi^{\prime}$\nobreakdash-\hspace{0pt}subgroup $\overline{K}$
of $\overline{G}$. Let $\overline{K}=B/N$, and let $R$ be a Hall
$\pi^{\prime}$\nobreakdash-\hspace{0pt}subgroup of $B$. It exists
because $B$ is a $\pi$\nobreakdash-\hspace{0pt}solvable group. Then $R$ is a
Hall $\pi^{\prime}$\nobreakdash-\hspace{0pt}subgroup of $G$, 
$RN/N=B/N= \overline{K}$ and by Lemma 1,
$$
N_{G}(R)N/N \subseteq N_{G/N}(RN/N)=N_{\overline{G}}(\overline{K}).
$$
By induction, $N_{G}(R)N/N \subseteq \overline{M}=M/N$, so $N_{G}(R)
\subseteq M$.  The Theorem 3 is proved.

Note that Theorem 2 follows from Theorem 3 if $\pi=\{p\}$. If $G$ is
a solvable group, then in the condition of Theorem 3,
$\pi(F(M/\mbox{Core}_{G}M))\ne \emptyset$, so
Theorem 1 follows from Theorem 3.

{\sc Proof of Corollary 1.1.} By Theorem 1, $G$ has a Sylow
$q$\nobreakdash-\hspace{0pt}subgroup $Q$ such that $N_{G}(Q)
\subseteq H$. Let $H_{1}$ be any subgroup such that $Q\subseteq
H_{1} \subseteq H$, and let $T=N_{G}(H_{1})$. By the Frattini lemma,
we have
$$
T=N_{T}(Q)H_{1} \subseteq N_{G}(Q)H_{1} \subseteq H, \
T=N_{G}(H_{1}) \subseteq H.
$$
The Corollary is proved.

\bigskip

{\bf References}

[1] Monakhov V.\,S.  Introduction to the Theory of Finite Groups and
Their Classes. Minsk: High school. 2006. (In Russian)

[2] Vedernikov V.\,A. About $\pi$-properties of a finite groups //
In the book: Arithmetic and subgroups a structure of finite
groups. Minsk: Science and technics. 1986. P. 13--19. (In Russian)

\bigskip

%\bigskip

\noindent D.\,V. GRITSUK

\noindent  Department of mathematics, Gomel F. Scorina State
University, Gomel 246019, BELARUS

\noindent E-mail address: Dmitry.Gritsuk@gmail.com

\bigskip

\noindent V.\,S. MONAKHOV

\noindent Department of mathematics, Gomel F. Scorina State
University, Gomel 246019,  BELARUS

\noindent E-mail address: Victor.Monakhov@gmail.com

\end{document}